\documentclass{amsart}
\usepackage{latexsym,amsmath,amssymb,graphics,amscd}
\@addtoreset{figure}{section}
\def\thefigure{\thesection.\@arabic\c@figure}
\def\fps@figure{h,t}
\@addtoreset{table}{bsection}
\def\thetable{\thesection.\@arabic\c@table}
\def\fps@table{h, t}
\@addtoreset{equation}{section}

\makeatother

\newcommand{\bfi}{\bfseries\itshape}
\newtheorem{theorem}{Theorem}[section]

\newtheorem{definition}[theorem]{Definition}

\newtheorem{lemma}[theorem]{Lemma}

\newtheorem{proposition}[theorem]{Proposition}
\newtheorem{remark}[theorem]{Remark}

\begin{document}
\title{On the symmetry breaking phenomenon}
\author{Petre Birtea, Mircea Puta, Tudor S. Ratiu and R\u{a}zvan Micu Tudoran}
\maketitle
\begin{abstract}
We investigate the problem of symmetry breaking in the framework
of dynamical systems with symmetry on a smooth manifold. Two cases
will be analyzed: general and Hamiltonian dynamical systems. We
give sufficient conditions for symmetry breaking in both cases.
\end{abstract}
\section{Introduction}
In this paper we investigate the problem of symmetry breaking in
the framework of dynamical systems with symmetry on a smooth
manifold. Two cases will be analyzed: general and Hamiltonian
dynamical systems. For general vector fields we translate the
problem locally to a symmetry breaking problem on a vector space
where the techniques are more developed and easier to apply
because of the existence of a Hilbert basis of invariant
polynomials for representations of compact Lie groups. The
symmetric steady state and Hopf bifurcations are considered as
examples.  Sufficient conditions for the existence of symmetry
breaking bifurcations are known in this  case. In the Hamiltonian
context we find sufficient conditions for the existence of
symmetry breaking bifurcations for Hamiltonian vector fields on
symplectic  vector spaces. The Hamiltonian steady state (passing
of eigenvalues of the linearization on the imaginary axis at the
origin) and the Hamiltonian Hopf bifurcations are analyzed.  It
turns out that the same methods as in the general case work,
because in the Hamiltonian context one can work only with the
invariant polynomials of the action and ignore the equivariant
polynomials. In addition, it  is shown how to locally translate
these Hamiltonian bifurcation problems from the manifold to the
tangent space at the bifurcation point. This is done by using the
equivariant Darboux theorem.

\section{Some basic results from the theory of group actions}
We shall need a few fundamental results form the theory of group
actions which we now briefly recall. For proofs and further
information see \cite{br}, \cite{dk}, \cite{kawakubo}, \cite{or}
and for the subsection on compact group representations see
\cite{k}, \cite{kc}, \cite{poenaru}.
\subsection{Twisted products}
Let $G$ be a Lie group and $H\subset G$ a Lie subgroup. Supose
that $H$ acts on the left on a manifold $A$. The {\bfi twisted
action\/} of $H$ on the product $G\times A$ is defined by
\[
h\cdot (g,a)=(gh,h^{-1}\cdot a), \quad h\in H, \quad g\in G, \quad
a\in A.
\]
Note that this action is free and proper by the freeness and
properness of the action on the $G$-factor. The {\bfi twisted
product\/} $G\times _{H}A$ is defined as the orbit space $(G\times
A)/H$ of the twisted action. The elements of $G\times _{H}A$ will
be denoted by $[g,a],$ $g\in G,$ $a\in A$. The twisted product
$G\times _{H}A$ is a $G$-space relative to the left action defined
by $g^{\prime }\cdot \lbrack g,a]=[g^{\prime }g,a]$. Also, the
action of $H$ on $A$ is proper if and only if the $G$-action on
$G\times _{H}A$ is proper. The isotropy subgroups of the
$G$-action on the twisted product $G\times _{H}A$ satisfy
\[
G_{[g,a]}=gH_{a}g^{-1}, \quad g\in G, \quad a\in A.
\]
\subsection{Slices}
Throughout this paragraph it will be assumed that $\Psi : G \times
Q \rightarrow Q $ is a left proper action of the Lie group $G$ on
the manifold $Q$. This action will not be assumed to be free, in
general. For $q\in Q$ we will denote by $H:=G_{q} := \{ g \in G
\mid g \cdot q = q \}$ the isotropy subgroup of the action $\Psi$
at $q $. We shall introduce also the following convenient
notation: if $K \subset G $ is a Lie subgroup of $G $ (possibly
equal to $G $), $\mathfrak{k}$ its Lie algebra,  and $q\in Q$,
then $\mathfrak{k}\cdot q := \{ \eta_Q(q) \mid \eta \in
\mathfrak{k} \}$ is the tangent space to the orbit $K\cdot q $ at
$q $. A {\bfi tube\/} around the orbit $G\cdot q$ is a
$G$-equivariant diffeomorphism $\varphi :G\times _{H}A\rightarrow
U$, where $U$ is a $G$-invariant neighborhood of $G\cdot q$ and
$A$ is some manifold on which $H$ acts. Note that the $G$-action
on the twisted product $ G\times _{H}A$ is proper since the
isotropy subgroup $H$ is compact and, consequently, its action on
$A$ is proper. Hence the $G$-action on $G\times _{H}A$ is proper.
Let $S$ be a submanifold of $Q$ such that $q\in S$ and $H\cdot
S=S$. We say that $S$ is a {\bfi slice\/} at $q$ if the map
\[
\varphi : G\times _{H}S \ni [ g,s] \mapsto g\cdot s \in  U
\]
is a tube about $G\cdot q$, for some $G$-invariant open
neighborhood of $G\cdot q$. Notice that if $S$ is a slice at $q$
then $g\cdot S$ is a slice at the point $g\cdot q$. The following
statements are equivalent:
\begin{enumerate}
\item[\textbf{(i)}] There is a tube $\varphi :G\times
_{H}A\rightarrow U$ about $G\cdot q$ such that $\varphi
([e,A])=S$. \item[\textbf{(ii)}] $S$ is a slice at $q$.
\item[\textbf{(iii)}] The submanifold $S$ satisfies the following
properties:
\begin{enumerate}
\item[\textbf{(a)}] The set $G\cdot S$ is an open neighborhood of
the orbit $G\cdot q$ and $S$ is closed in $G\cdot S$.
\item[\textbf{(b)}] For any $s\in S$ we have
$T_{s}Q=\mathfrak{g}\cdot s+T_{s}S$. Moreover, $\mathfrak{g}\cdot
s\cap T_{s}S=\mathfrak{h}\cdot s$, where $\mathfrak{h}$ is the Lie
algebra of $H $. In particular $T_{q}Q=\mathfrak{g}\cdot q\oplus
T_{q}S$. \item[\textbf{(c)}] $S$ is $H$-invariant. Moreover, if
$s\in S$ and $g\in G$ are such that $g\cdot s\in S$, then $g\in
H$. \item[\textbf{(d)}] Let $\sigma :U\subset G/H\rightarrow G$ be
a local section of the submersion $G\rightarrow G/H$. Then the map
$F:U\times S\rightarrow Q$ given by $F(u,s):=\sigma (u)\cdot s$ is
a diffeomorphism onto an open set of $Q$.
\end{enumerate}
\item[\textbf{(iv)}] $G\cdot S$ is an open neighborhood of $G\cdot
q$ and there is an equivariant smooth retraction $r:G\cdot
S\rightarrow G\cdot q$ of the injection $G\cdot q\hookrightarrow
G\cdot S$ such that $r^{-1}(q)=S$.
\end{enumerate}
\begin{theorem}
\textbf{(Slice Theorem)} Let $Q$ be a manifold and $G$ be a Lie
group acting properly on $Q$ at the point $q\in Q$. Then there
exists a slice for the $G$-action at $q$.
\end{theorem}
\begin{theorem}
\textbf{(Tube Theorem)} Let $Q$ be a manifold and $G$ a Lie group
acting properly on $Q$ at the point $q\in Q$, $H:=G_{q}$. Then
there is a tube $ \varphi :G\times _{H}B\rightarrow U$ about
$G\cdot q$ such that $\varphi ([e,0])=q$, $\varphi ([e,B])=:S$ is
a slice at $q$; $B$ is an open  $H$-invariant neighborhood of $0$
in the vector space $T_{q}Q/T_{q}(G\cdot q) $, on which $H$ acts
linearly by $h\cdot (v_{q}+T_{q}(G\cdot q)):=T_{q}\Psi
_{h}(v_{q})+T_{q}(G\cdot q)$.
\end{theorem}
If $Q$ is a Riemannian manifold then $B$ can be chosen to be a
$G_{q}$-invariant neighborhood of $0$ in $(\mathfrak{g}\cdot
q)^{\perp }$, the orthogonal complement to $\mathfrak{g}\cdot q$
in $T_{q}Q$. In this case $ U=G\cdot
\operatorname{\operatorname{Exp}}_{q}(B)$, where
$\operatorname{\operatorname{Exp}}_{q}: T_q Q \rightarrow Q $ is
the Riemannian exponential map.
\subsection{Type submanifolds and fixed point subspaces}
Let $G$ be a Lie group acting on a manifold $Q$. Let $H$ be a
closed subgroup of $G$. We define the following subsets of $Q$ :
\begin{eqnarray*}
Q_{(H)} &=&\{q\in Q\mid G_{q}=gHg^{-1},g\in G\}, \\
Q^{H} &=&\{q\in Q\mid H\subset G_{q}\}, \\
Q_{H} &=&\{q\in Q\mid H=G_{q}\}.
\end{eqnarray*}
All these sets are submanifolds of $Q $. The set $Q_{(H)}$ is
called the $(H)$-{\bfi orbit type submanifold\/}, $ Q_{H}$ is the
$H$-{\bfi isotropy type submanifold\/}, and $Q^{H}$ is the
$H$-{\bfi fixed point submanifold\/}. We will collectively call
these subsets the {\bfi type submanifolds}. We have:
\begin{itemize}
\item  $Q^{H}$ is closed in $Q$; \item $Q_{(H)}=G\cdot Q_{H}$;
\item $Q_H $ is open in $Q^H $. \item the tangent space at $q \in
Q_H $ to $Q_H $ equals
\[
T_{q}Q_{H}=\{v_{q}\in T_{q}Q\mid T_{q}\Psi_{h}(v_{q}) =v_{q},\,
\text{~for~all~} h\in H\}=(T_q Q)^H=T_q Q^H;
\]
\item $T_{q}(G\cdot q)\cap (T_{q}Q)^{H}=T_{q}(N(H)\cdot q)$, where
$N(H) $ is the normalizer of $H $ in $G $; \item if $H$ is compact
then $Q_{H}=Q^{H}\cap Q_{(H)}$ and $Q_H $ is closed in $Q_{(H)}$.
\end{itemize}
If $Q$ is a vector space on which $H$ acts linearly, the set
$Q^{H}$ is called in the physics literature  the {\bfi space of
singlets\/} or the {\bfi space of invariant vectors\/}.
\begin{theorem}
\label{Teorema de stratificare} \textbf{(The stratification
theorm)} Let $Q$ be a smooth manifold and $G$ be a Lie group
acting properly on it. The connected components of the orbit type
manifolds $Q_{(H)}$ and their projections onto the orbit space
$Q_{(H)}/G$ constitute a Whitney stratification of $Q$ and $Q/G$,
respectively. This stratification of $Q/G$ is minimal among all
Whitney stratifications of $Q/G$.
\end{theorem}
The proof of this result, that can be found in \cite{dk} or
\cite{pflaum}, is based on the Slice Theorem and on a series of
extremely important properties of the orbit type manifolds
decomposition that we enumerate in what follows. We start by
recalling that the set of conjugacy classes of subgroups of a Lie
group $G$ admits a partial order by defining $ (K)\preceq (H)$ if
and only if $H$ is conjugate to a subgroup of $K$. Also, a point
$q\in Q$ in a proper $G$-subspace of $Q$ (or its corresponding
$G$-orbit, $ G\cdot q$) is called {\bfi principal\/} if its
corresponding local orbit type manifold is open in $Q$. The orbit
$G\cdot q$ is called {\bfi regular\/} if the dimension of the
orbits nearby coincides with the dimension of $ G\cdot q$. The set
of principal and regular orbits will be denoted by $ Q_{princ}/G$
and $Q_{reg}/G$, respectively. Using this notation we have:
\begin{itemize}
\item  For any $q\in Q$ there exists a neighborhood $U$ of $q$
that intersects only finitely many connected components of
finitely many orbit type manifolds. If $Q$ is compact or a linear
space where $G$ acts linearly, then the $G$-action on $Q$ has only
finitely many distinct connected components of orbit type
manifolds.

\item  For any $q\in Q$ there exists an open neigborhood $U$ of
$q$ such that $(G_{q})\preceq (G_{x})$, for all $x\in U$. In
particular, this implies that $\dim G\cdot q\leq \dim G\cdot x$,
for all $x\in U$.

\item  {\bfi Principal Orbit Theorem}: For every connected
component $ Q^{0}$ of $Q$ the subset $Q_{princ}\cap Q^{0}$ is
connected, open, and dense in $Q^{0}$. Each connected component
$(Q/G)^{0}$ of $Q/G$ contains only one principal orbit type, which
is connected open and dense in $(Q/G)^0 $.
\end{itemize}
\subsection{General facts about compact Lie groups representations}
Let $G$ be a compact Lie group acting on a finite dimensional
vector space $V$. The compactness of $G$ implies the existence of
a $G$-invariant inner product $\langle \!\langle \cdot, \cdot
\rangle\!\rangle $ on $V$, so one can always assume that the
$G$-representation is orthogonal. Denote by $\mathcal{P}_{G}(V)$
the ring of $G$-invariant real valued polynomials on $V$. The
Hilbert-Weyl theorem implies that $\mathcal{P}_{G}(V)$ is finitely
generated over $\mathbb{R}$ and it can be shown that one can
choose the generators to be homogeneous polynomials. The elements
$\{\theta_{1},...,\theta_{l}\}$ of such a basis over $\mathbb{R}$
are called {\bfi invariant generators}. The following theorem of
Schwarz gives a method of writing a $G$-equivariant smooth map in
terms of invariants and equivariants.
\begin{theorem}
Let $V$ and $W$ be two $G$-representation spaces. The set of
$G$-equivariant polynomial maps from $V$ to $W$ is a finitely
generated $\mathcal{P}_{G}(V)$-module. Let $\{F_{1},...,F_{k}\}$
be one of its minimal bases. Then, if $F$ is a smooth
$G$-equivariant map from $V$ to $W$, there exist smooth functions
$f_{1},...,f_{k}$ on $\mathbb{R}^{l}$ such that
\begin{equation*}
F(v)=\sum_{i=1}^{k}f_{i}(\theta_{1}(v),...,\theta_{l}(v))F_{i}(v).
\end{equation*}
\end{theorem}
The next theorem gives an identification between the orbit space
of a linear action of a compact Lie group $G$ on a vector space
$V$ and a semialgebraic subset of $\mathbb{R}^{l}$.
\begin{theorem}
Let $\pi$ be the Hilbert map $V \ni v
\mapsto{(\theta_{1}(v),...,\theta_{l}(v))}\in{\mathbb{R}^{l}}$.
Then $\pi$ induces a homeomorphism between the orbit space $V/G$
and $\operatorname{Im}(\pi)$.
\end{theorem}
The invariant generators are intimately related to the orthogonal
space to the $G $-orbits. If $v \in V $, let $N_v : = \{u \in V
\mid \langle \!\langle u, \xi \cdot v \rangle\!\rangle = 0,
\text{for all}\; \xi\in \mathfrak{g}\}$ be the orthogonal space to
the orbit $G \cdot v $. Define the {\bfi normal invariant
subspace} $N_v^{(0)} : = \{u \in N_v \mid g \cdot u = u, \text{for
all} \; g \in G_v\}\subset N_v$. The following theorem is proved
in \cite{as}.
\begin{theorem}
\label{as_theorem} The vector subspace of $V$ generated by
$\{\operatorname{grad} \theta_i (v) \mid i = 1, \dots, l\}$
coincides with the normal invariant subspace $N_v^{(0)}$.
\end{theorem}
Two useful consequences of this fact are the orthogonal sum
decomposition
\[
T_v V_{(H)} = T_v(G\cdot v) \oplus N_v^{(0)}
\]
and  the fact that $F_i(v) \in T_v\left(V_{(G_v)}\right) $ for any
equivariant generator $F_i $, $i = 1, \dots , k $ of the
$G$-representation.
\medskip
Since the orbit space is the image of a finite dimensional vector
space by a polynomial map, the Tarski-Seidenberg theorem
guarantees that it is a semialgebraic set. Thus, as a
semialgebraic set, $V/G$ has a canonical stratification. A theorem
of Bierstone implies that this canonical semialgebraic
stratification is the same as the stratification by orbit types
obtained by projection of those from $V$. Since the orbit space
$V/G $ inherits a smooth structure from $V$ and $G$, but in
general is not a manifold, we have to define the meaning of a
smooth vector field on $V/G $.
\begin{definition}
Let $X$ be a smooth vector field defined on a neighborhood $U$ of
$0\in\mathbb{R}^{l}$ and tangent to the strata of $\pi(V)\cap{U}$.
Then $X$ is said to be {\bfi smooth vector field on\/} $V/G$ in a
neighborhood of the origin.
\end{definition}
If one projects a smooth $G$-equivariant vector field on the orbit
space, one obtains a smooth vector field on $V/G$ in the sense
defined above. When we speak of a vector field on $V/G$ we mean a
vector field which is tangent to each stratum of $V/G$. The next
step gives a method for projecting a $G$-equivariant vector field
on the orbit space. For doing this, let $F$ be a smooth
$G$-equivariant vector field on the finite dimensional vector
space $V$. Let $\{\theta_{1},...,\theta_{l}\}$ be a set of
generators for $\mathcal{P}_{G}(V)$ and let $\{F_{1},...,F_{k}\}$
be a set of generators for the $\mathcal{P}_{G}(V)$-module of
$G$-equivariant polynomials maps from $V$ to $V$. By the Schwarz
theorem we can write $F$ in terms of invariants and equivariants
as
\begin{equation*}
F(v)=\sum_{i=1}^{k}f_{i}(\theta_{1}(v),...,\theta_{l}(v))F_{i}(v),
\end{equation*}
where the functions $f_{i}$ are smooth functions on
$\mathbb{R}^{l}$. To obtain the projection of this vector field on
the orbit space $V/G$, one has to compute the derivatives of the
invariant generators $\theta_{i}$ with respect to time:
\begin{align*}
\dot{\theta}_{j}&
=\langle\!\langle{\operatorname{grad}_{v}(\theta_{j})},\dot{v}\rangle
\!\rangle =\langle\!\langle{\operatorname{grad}_{v}(\theta_{j})},
\sum_{i=1}^{k}f_{i}(\theta(v))F_{i}(v)\rangle\!\rangle
\\ &=\sum_{i=1}^{k}f_{i}(\theta(v))\langle\!\langle{\operatorname{grad}_{v}
(\theta_{j})},F_{i}(v)\rangle\!\rangle.
\end{align*}
Since the expression
$\langle\!\langle{\operatorname{grad}_{v}(\theta_{j})},F_{i}(v)\rangle\!\rangle$
is $G$-invariant one can rewrite it in terms of the invariant
generators. Thus the projection $\widetilde{F}$ of $F$ on the
orbit space $V/G$ has the expression
\begin{equation}
\widetilde{F}(\theta)=\sum_{i=1}^{k}f_{i}(\theta){\widetilde{F}}_{i}(\theta),
\end{equation}
where if $v$ is such that $\pi(v)=\theta$ and
\begin{equation*}
{\widetilde{F}}_{i}(\theta):=
\begin{bmatrix}
\langle\!\langle{F_{i}(v),\operatorname{grad}_{v}\theta_{1}(v)}\rangle\!\rangle \\
\vdots\\
\langle\!\langle{F_{i}(v),\operatorname{grad}_{v}\theta_{l}(v)}\rangle\!\rangle
\end{bmatrix}
\quad \text{~for~} \quad i = 1, \dots, k.
\end{equation*}

\section{Formulation of Symmetric Systems on Vector Spaces}
This section presents a method of translating a bifurcation
problem for a general vector fields family with $G$-symmetry from
a smooth manifold to a vector space. We begin by quickly recalling
the necessary notations and statements for a left smooth proper
Lie group action $\Phi: G \times Q \rightarrow Q $, where $Q $ is
a smooth manifold and $G $ is a Lie group with Lie algebra
$\mathfrak{g}$.
\begin{definition}
A vector field $X:Q\rightarrow TQ$ is said to be $G$-equivariant
if
\[
T_{q}\Phi _{g}(X(q))=X(\Phi _{g}(q))
\]
for all $q\in Q$ and $g\in G$. If $X$ is $G$-equivariant, then $G$
is said to be a \textbf{symmetry group} of the dynamical system
$\dot{q}=X(q)$.
\end{definition}
\begin{definition}
A \textbf{relative equilibrium} of a $G$-equivariant vector field
$X$ is a point $q_{e}\in Q$ at which the values of $X$ and of the
infinitesimal generator $\xi_Q$ of some element $\xi \in
\mathfrak{g}$, called {\bfi velocity\/} of the relative
equilibrium,  coincide, that is,
\[
X(q_{e})=\xi _{Q}(q_{e}) : = \left.\frac{d}{dt}\right|_{t=0}
\Phi(\exp t \xi, q_e),
\]
where $\exp: \mathfrak{g} \rightarrow G $ is the Lie group
exponential map. A relative equilibrium $q_{e}$ is said to be
\textbf{asymmetric }if $ \mathfrak{g}_{q_{e}}=\{0\}$, and
\textbf{symmetric} otherwise.
\end{definition}
From the definition above and the uniqness of the integral curves
of a vector field with given initial condition, it follows that
the dynamical orbit $q(t) $ of $X $ staring at the relative
equilibrium $q_{e}$, that is, $q(0) = q_e $,  coincides with the
curve $\exp t \xi \cdot q _e : = \Phi(\exp t \xi, q_e)$ given by
the group action. Note that if $\xi$ is a velocity of the relative
equilibrium $q_e $, then so is $\xi+ \eta$, where $\eta \in
\mathfrak{g}_{q_e}:=
\{\zeta\in{\mathfrak{g}}\mid{\zeta_{Q}}(q_{e})=0_{q_{e}}\}$. The
Lie subalgebra $\mathfrak{g}_{q_e}$ is called the {\bfi symmetry
algebra\/} of $q_e $. Denote for each $g \in G $ by
$\operatorname{Ad}: G \times \mathfrak{g}\rightarrow \mathfrak{g}$
the adjoint representation of $G $ on $\mathfrak{g}$; for each $g
\in G $, $\operatorname{Ad}_g$ is the Lie algebra isomorphism
obtained by taking the derivative of the conjugation isomorphism
on $G $.
\begin{proposition}
If $q_{e}$ is a relative equilibrium with velocity $\xi \in
\mathfrak{g}$, then, for any $g\in G$, $\Phi_{g}(q_{e})$ is also a
relative equilibrium with velocity $\operatorname{Ad}_{g}\xi $.
\end{proposition}
\noindent\textbf{Proof.} Indeed, since $X $ is $G $-equivariant,
we get
\begin{equation*}
X(\Phi_{g}(q_{e}))=T_{q_{e}}\Phi_{g}(X(q_{e}))
=\left(T_{q_{e}}\Phi_{g} \circ \xi_{Q}\circ
\Phi_{g^{-1}}\right)\left(\Phi_{g}(q_{e})\right)
=(\operatorname{Ad}_{g}\xi)_{Q}(\Phi_{g}(q_{e})),
\end{equation*}
which proves the statement.\quad $\blacksquare$

By one of the remarks following Theorem \ref{Teorema de
stratificare}, we can choose a $G$-invariant tubular open
neigborhood $W_{q_{e}}$ of the orbit $G\cdot{q_{e}}$ such that the
isotropy subgroup $G_x $ for any $x \in W_{q_e} $ is conjugate to
a (possibly non-strict) subgroup of $G_{q_{e}}$. The next lemma
will show that we can do better.
\begin{lemma}
\label{ONB} Let $L$ be a compact Lie group acting on a smooth
manifold $Q$. Assume that $q\in Q$ is a fixed point of the
$L$-action. Then any open neighborhood of $q$ contains a
$L$-invariant open neighborhood of $q$.
\end{lemma}

\noindent\textbf{Proof.} Let $\Phi:L\times{Q}\to{Q}$ be the group
action and $U$ an arbitrary neighborhood of $q$. The set
$\Phi^{-1}(U)$ is clearly open and contains $L\times\{q\}$. By
continuity of $\Phi$, for any $g\in {L}$, there are open
neighborhoods $W_{g}$ of $g$ in $L$ and $V_{g}$ of $q$ in $Q$ such
that $W_{g}\times{V_{g}}\subset\Phi^{-1}(U)$. Since, by
hypothesis, $L$ is compact, the open cover $\{W_{g }\mid g\in
{L}\}$ of $L$ has a finite subcover $\{W_{g_{1}},...,W_{g_{n}}\}$.
Let $V:=\bigcap_{i=1}^{n}V_{g_{i}}$. The $L$-invariant set
$W:=\Phi(L,V)\subset U$ is clearly open since $W=
\bigcup_{g\in{L}}{\Phi_{g}(V)}$ and, by construction, contains the
point $q$. \quad $\blacksquare$

Properness of the $G $-action on $Q$ implies that the isotropy
subgroup $G_{q_{e}}$ is compact. By the lemma above we can
conclude that there exists an open $G_{q_{e}}$-invariant
neighborhood $W_{q_{e}}^{1}$ of $q_{e}$ contained in $W_{q_{e}}$.
Thus the local study of the symmetry breaking from $G_{q_{e}}$ to
a smaller symmetry group can be done by restricting the action to
$W_{q_{e}}^{1}$. Next, assume that $Q $ carries a $G $-invariant
Riemannian metric. The Riemannian exponential map
$\operatorname{Exp}: TQ\rightarrow Q $ is also $G $-equivariant,
relative to the naturally lifted $G$-action on $TQ$. Since
$\operatorname{Exp}_{q_e}:T_{q_e}Q \rightarrow Q $ is a local
diffeomorphism around the origin, we can further shrink the
neighborhood $W_{q_{e}}^{1}$ (keeping the $G_{q_e}$-invariance),
till it is contained in the image under $\operatorname{Exp}_{q_e}$
of the neighborhood of the origin where $\operatorname{Exp}_{q_e}$
is a diffeomorphism. Denote this shrunk $G_{q_e}$-invariant
neighborhood by $U_{q_e} $. The inverse map of the Riemannian
exponential on $U_{q_{e}}$, $\operatorname{Log}_{q_{e}}: U_{q_{e}}
\to \operatorname{Log}_{q_e}(U_{q_{e}}) \subset T_{q_{e}}Q $ is
therefore also $G_{q_{e}}$-equivariant and so its image
$\operatorname{Log}_{q_{e}}(U_{q_{e}})=:V_{0_{q_{e}}}$ will be a
$G_{q_{e}}$-invariant open neighborhood of
$0_{q_{e}}\in{T_{q_{e}}Q}$.

\begin{proposition}
Let $(Q,\langle\!\langle \cdot, \cdot \rangle\!\rangle)$ be a
smooth Riemannian manifold, $q_{e}\in{Q}$, $G$ a Lie group acting
properly and smoothly on $Q$ by isometries,  and
$X\in\mathfrak{X}(Q)$ a smooth $G$-equivariant vector field. Then,
in the notations above, the vector field
$\operatorname{Exp}_{q_{e}}^{*}\left(X|_{U_{q_e}}\right)
\in{\mathfrak{X}(V_{0_{q_{e}}})}$ is a smooth
$G_{q_{e}}$-equivariant vector field on the open neighborhood
$V_{0_{q_e}}$ of the origin $0_{q_{e}}$ in the vector space
$T_{q_{e}}Q$.
\end{proposition}

\noindent\textbf{Proof.} The proof is a simple consequence of the
fact that $\operatorname{Exp}_{q_{e}}: V_{q_e} \subset T_{q_{e}}Q
\rightarrow U_{q_e} \subset Q$ and $\operatorname{Log}_{q_e}:
U_{q_e} \rightarrow V_{q_e}$ are $G_{q_{e}}$-equivariant
diffeomorphisms inverse to each other. Indeed, since $X$ is
$G$-equivariant, its restriction $X|_{U_{q_e}}$ is
$G_{q_{e}}$-equivariant. Therefore, for any $h \in G_{q_{e}}$ and
$v_{q_e} \in V_{q_e}$ we have
\begin{align*}
&\left(\operatorname{Exp}_{q_{e}}^{*}\left(X|_{U_{q_e}}\right)\right)(T_{q_e}
\Phi_{h}(v_{q_e})) \\
&\qquad =T_{{\operatorname{Exp}_{q_{e}}}(T_{q_e}{\Phi_{h}}
(v_{q_e}))}\operatorname{Log}_{q_{e}}
\left({X(\operatorname{Exp}_{q_e}(T_{q_e}{\Phi_{h}}(v_{q_e})))}\right)\\
&\qquad=T_{\Phi_{h}({\operatorname{Exp}_{q_{e}}}(v_{q_{e}}))}\operatorname{Log}_{q_{e}}
\left({X(\Phi_{h}({\operatorname{Exp}_{q_{e}}(v_{q_{e}})}))}\right)\\
&\qquad=T_{\operatorname{Exp}_{q_{e}}(v_{q_{e}})}
(\operatorname{Log}_{q_{e}}\circ{\Phi_{h}})\left({X(\operatorname{Exp}_{q_{e}}
(v_{q_{e}}))}\right)\\
&\qquad=T_{\operatorname{Exp}_{q_{e}}(v_{q_{e}})}({T_{q_{e}}\Phi_{h}}
\circ{\operatorname{Log}_{q_{e}}})\left({X(\operatorname{Exp}_{q_{e}}(v_{q_{e}}))}\right)\\
&\qquad=T_{v_{q_{e}}}(T_{q_{e}}\Phi_{h})\left(\left(\operatorname{Exp}_{q_{e}}^{*}
\left(X|_{U_{q_e}}\right)\right)(v_{q_{e}})\right),
\end{align*}
which proves that
$\operatorname{Exp}_{q_{e}}^{*}\left(X|_{U_{q_e}}\right)$ is a
$G_{q_{e}}$-equivariant vector field on $V_{0_{q_{e}}}$.
\quad$\blacksquare$

\begin{proposition}
\label{correspondence statement} Let $(Q,\langle\!\langle \cdot,
\cdot \rangle\!\rangle)$ be a smooth Riemannian manifold, $G$ a
Lie group acting properly on $Q$ by isometries, and
$X_{\lambda}\in\mathfrak{X}(Q)$, $\lambda\in{\mathbb{R}}$, a
family of smooth $G$-equivariant vector fields. Suppose that
$q_{e}\in Q$ is a (relative) equilibrium for $X_{\lambda_{0}}$
with symmetry $G_{q_{e}}$. Assume the family of vector fields
$\operatorname{Exp}_{q_{e}}^{*}X_{\lambda} \in
\mathfrak{X}(V_{0_{q_e}})$ has the property that there exist
$\epsilon>0$ and a continuous curve
$\gamma:[\lambda_{0},\lambda_{0}+\epsilon]\to{V_{0_{q_{e}}}}$,
$\gamma(\lambda_0) = 0_{q_e}$, such that all the points on the
curve $\gamma(\lambda) $ for $\lambda>\lambda_{0}$ are relative
equilibria of $\operatorname{Exp}_{q_{e}}^{*}X_{\lambda}$ with
isotropy subgroups conjugate to $H \subset G_{q_{e}}$.  Then the
curve
$\operatorname{Exp}_{q_{e}}\circ\gamma:[\lambda_{0},\lambda_{0}+\epsilon]\to{U_{q_{e}}}$
is formed by relative equilibria of $X_{\lambda}$ with isotropy
subgroups conjugate to $H$ for all $\lambda>\lambda_{0}$ and to
$G_{q_e}$ for $\lambda=\lambda_{0}$.
\end{proposition}

\noindent\textbf{Proof.} Let
$X_{\lambda}\in\mathfrak{X}(U_{q_{e}})$ be one of the vector
fields restricted to the $G_{q_{e}}$-invariant open neighborhood
$U_{q_{e}}$. Suppose that $\gamma(\lambda)$ is a relative
equilibrium of the vector field
$\operatorname{Exp}_{q_{e}}^{*}X_{\lambda}$. This means that there
exists an element
$\xi(\lambda)\in{\mathfrak{g}_{q_{e}}}\subset{\mathfrak{g}}$ such
that
\begin{align*}
(\operatorname{Exp}_{q_{e}}^{*}X_{\lambda})
(\gamma(\lambda))=(\xi(\lambda))_{T_{q_{e}}Q}(\gamma(\lambda)),
\end{align*}
that is,
\begin{align*}
T_{\operatorname{Exp}_{q_{e}}(\gamma(\lambda))}
\operatorname{Log}_{q_{e}}\left({X_{\lambda}(\operatorname{Exp}_{q_{e}}
(\gamma(\lambda)))}\right)=
(\xi(\lambda))_{T_{q_{e}}Q}(\gamma(\lambda)).
\end{align*}
If we apply to the above equality
$T_{\gamma(\lambda)}\operatorname{Exp}_{q_{e}}$ and use the
$G_{q_{e}}$-equivariance of the exponential map we get
\begin{align*}
X_{\lambda}\left(\operatorname{Exp}_{q_{e}}(\gamma(\lambda))\right)
&= T_{\gamma(\lambda)}\operatorname{Exp}_{q_{e}}
\left((\xi(\lambda))_{T_{q_{e}}Q}(\gamma(\lambda))\right) \\&=
(\xi(\lambda))_{Q}(\operatorname{Exp}_{q_{e}}(\gamma(\lambda))).
\end{align*}
This means that $\operatorname{Exp}_{q_{e}}(\gamma(\lambda))$ is a
relative equilibrium of $X_{\lambda}$. So we have proved that the
curve $\operatorname{Exp}_{q_{e}}\circ{\gamma}$ is formed by
relative equilibria of $X_\lambda$ for any
$\lambda\in[\lambda_{0},\lambda_{0}+\epsilon]$. Because the
exponential map is a $G_{q_{e}}$-equivariant diffeomorphism on the
neighborhoods considered, we have
$G_{\operatorname{Exp}_{q_{e}}(\gamma(\lambda))} =
G_{\gamma(\lambda)}$, and so all the points on the curve
$\gamma(\lambda) $ with the possible exception of
$\gamma(\lambda_{0})$ are conjugate to $H$. \quad $\blacksquare$

Since the $G_{q_e}$-equivariant family of vector fields
$\operatorname{Exp}_{q_{e}}^{*}X_{\lambda}\in\mathfrak{X}(V_{0_{q_{e}}})$
are defined on an open $G_{q_e}$-invariant neighborhood of the
origin in the vector space $T_{q_e}Q $, we can arbitrarily extend
their restrictions to a smaller $G_{q_e}$-invariant neighborhood
of the origin to smooth $G_{q_e}$-equivariant vector fields
$\overline{X_{\lambda}}\in\mathfrak{X}(T_{q_{e}}Q)$ on the whole
space $T_{q_{e}}Q$. However, on $T_{q_{e}}Q$ we have an orthogonal
representation of the compact group $G_{q_e}$ and thus, as
explained earlier, the vector fields $\overline{X_{\lambda}} \in
\mathfrak{X}(T_{q_{e}}Q)$ can be expressed in terms of the
invariant and equivariant polynomials on $T_{q_{e}}Q$. Thus to
study the local symmetry breaking phenomenon around the point $q_e
\in Q $, one can translate the problem to the vector space
$T_{q_{e}}Q$ and study this problem for the $G_{q_e}$-equivariant
family $\overline{X_{\lambda}} \in \mathfrak{X}(T_{q_e}Q)$.

\section{Stationary and Hopf Bifurcation with Symmetry}
We shall analyze now the local symmetry breaking phenomenon around
$0_{q_{e}}\in T_{q_{e}}Q$ for the family of $G_{q_e}$-equivariant
vector fields $\overline{X_{\lambda}} \in
\mathfrak{X}(T_{q_{e}}Q)$ under the following two nondegeneracy
conditions:
\begin{itemize}
\item \textit{Stationary bifurcation with symmetry}:
\begin{equation*}
D\overline{X_{\lambda}}(0_{q_{e}})={\sigma}_{X}(\lambda)\operatorname{Id}_{T_{q_{e}}Q}
\end{equation*}
where $\sigma_{X}:(-\epsilon, \epsilon) \rightarrow \mathbb{R}$ is
a smooth map satisfying  ${\sigma}_{X}(0)=0$ and
$\sigma^{\prime}_{X}(0)\neq 0$. In addition, it is assumed that
the compact group $G_{q_e}$ acts {\bfi absolutely irreducibly\/}
on $T_{q_e}Q $, i.e., that only scalar multiples of the identity
commute with the $G_{q_e}$-action. \item \textit{Hopf bifurcation
with symmetry}:
\begin{equation*}
D\overline{X_{\lambda}}(0_{q_{e}})={\sigma}_{X}(\lambda)
\operatorname{Id}_{T_{q_{e}}Q}+{\rho}_{X}(\lambda)
\mathbb{J}_{T_{q_{e}}Q}
\end{equation*}
where $\sigma_{X}, \rho_{X} :(-\epsilon, \epsilon) \rightarrow
\mathbb{R}$ are smooth maps satisfying  ${\sigma}_{X}(0)=0$,
$\sigma^{\prime}_{X}(0)\neq 0$, and ${\rho}_{X}(0)\neq 0$. In
addition, it is assumed that the vector space $T_{q_e}Q $ carries
a natural complex structure
$\operatorname{\mathbb{J}}_{T_{q_{e}}Q}$ and that the compact
group $G_{q_e}$ acts on this complex vector space $T_{q_e} Q$
irreducibly (as a complex representation). In this case the
manifold $Q $ is necessarily of even dimension.
\end{itemize}

In the following will denote by
$\mathcal{C}^{\infty}(T_{q_{e}}Q\times
\mathbb{R},T_{q_{e}}Q)^{G_{q_{e}}}_{o}$ the class of
$G_{q_{e}}$-equivariant vector fields that verifies one of the
above nondegeneracy conditions.
\medskip

Let $\{\theta_{1}(v_{q_{e}}),...,\theta_{l}(v_{q_{e}})\}$ be a set
of generators for $\mathcal{P}_{G_{q_{e}}}(T_{q_{e}}Q)$ and let
$\{F_{1},...,F_{k}\}$ be a set of generators for the
$\mathcal{P}_{G_{q_{e}}}(T_{q_{e}}Q)$-module of
$G_{q_{e}}$-equivariant maps from $T_{q_{e}}Q$ to itself. Then by
the Schwarz theorem we have
\begin{equation}
\label{x_bar}
\overline{X_{\lambda}}(v_{q_{e}})=\sum_{j=1}^{k}f_{j}(\theta_{1}(v_{q_{e}}),...,
\theta_{l}(v_{q_{e}}),\lambda)F_{j}(v_{q_{e}}),
\end{equation}
where the functions $f_{i}$ are smooth functions on
$\mathbb{R}^{l}$. The projection of this family of vector fields
on the orbit space $(T_{q_{e}}Q)/G_{q_{e}}$ has the expression:
\begin{equation}
\label{theta_dot_eq}
\dot{\theta}=\sum_{i=1}^{k}f_{i}(\theta,\lambda){\widetilde{F}}_{i}(\theta),
\end{equation}
where, if $v_{q_{e}} \in T_{q_e} Q$ is such that
$\pi(v_{q_{e}})=\theta \in ( T_{q_e} Q)/G_{q_e}$, we have
\begin{equation*}
{\widetilde{F}}_{i}(\theta)=
\begin{bmatrix}
\langle\!\langle{F_{i}(v_{q_{e}}),\operatorname{grad}_{v_{q_{e}}}
\theta_{1}(v_{q_{e}})}\rangle\!\rangle \\
\vdots\\
\langle\!\langle{F_{i}(v_{q_{e}}),\operatorname{grad}_{v_{q_{e}}}
\theta_{l}(v_{q_{e}})}\rangle\!\rangle
\end{bmatrix} \quad \text{for} \quad i = 1, \dots, k.
\end{equation*}
According to the orbit space projection of the vector field
$\overline{X_{\lambda}}$ we introduce the maps
\begin{equation}
\label{definition of tilde f} \widetilde{F}:(T_{q_{e}}Q)/G_{q_{e}}
\times \mathbb{R}^{k} \ni (\theta,t) \mapsto
\sum_{i=1}^{k}t_{i}\widetilde{F}_{i}(\theta) \in
(T_{q_{e}}Q)/G_{q_{e}}
\end{equation}
\begin{equation*}
\gamma:{\mathbb{R}} \ni \lambda
\mapsto{(f_{1}(0,\lambda),...,f_{k}(0,\lambda))}\in{\mathbb{R}^{k}}
\end{equation*}
and define the sets:
\begin{equation*}
\widetilde{\mathcal{E}}=\left\{(\theta,t)\in(T_{q_{e}}Q)/G_{q_{e}}\times{\mathbb{R}^{k}}\mid
{\widetilde{F}(\theta,t)=0}\right\}
\end{equation*}
\begin{equation*}
\widetilde{\mathcal{E}}_{(H)}=\left\{(\theta,t)\in\widetilde{\mathcal{E}}\mid
{\theta\in\Pi\left((T_{q_{e}}Q)_{(H)}\right)}\right\},
\end{equation*}
\begin{equation*}
\widetilde{\mathcal{A}}_{(H)}=\overline{\widetilde{\mathcal{E}}_{(H)}}
\cap\left(\left\{[0]\right\}\times{\mathbb{R}^{k}}\right),
\end{equation*}
where $\Pi:T_{q_{e}}Q\rightarrow(T_{q_{e}}Q)/G_{q_{e}}$ is the
orbit space projection and $H \subset G_{q_{e}}$ is an isotropy
subgroup of the $G_{q_{e}}$-action on $T_{q_{e}}Q$. Recall that
$\widetilde{\mathcal{E}}$ is a semialgebraic subset of
$\mathbb{R}^{l}\times{\mathbb{R}^{k}}$. Denote by $\Sigma$ its
canonical Whitney stratification. This stratification induces a
Whitney semialgebraic stratification $\Sigma_{(H)}$ on each
$\widetilde{\mathcal{E}}_{(H)}$. Let $\mathcal{B}$ denote the
Whitney stratification $\Sigma_{(G_{q_{e}})}$ of
$\widetilde{\mathcal{E}}_{(G_{q_{e}})}=\{0\}\times{\mathbb{R}^{k}}$.
Next we will define a notion of $G_{q_{e}}$-transversality
introduced by K\oe nig and Chossat \cite{kc} that was inspired by
Bierstone \cite{bierstone} and Field \cite{field}.
\begin{definition}
\label{transversality definition} The projection to the orbit
space $(T_{q_{e}}Q)/G_{q_{e}}$ of a $G_{q_{e}}$-equivariant smooth
family of vector fields $\overline{X_{\lambda}} \in
\mathfrak{X}(T_{q_e} Q)$  is {\bfi transverse\/} to
$[0_{q_e}]\in{(T_{q_{e}}Q)/G_{q_{e}}}$ at $\lambda=
0\in{\mathbb{R}}$ if $\gamma: \mathbb{R} \rightarrow \mathbb{R}^k$
is transverse at $\lambda= 0 $ to the stratum of $\mathcal{B}$
containing $\gamma(0) $, that is, $\operatorname{im}T_0 \gamma +
T_{\gamma(0)} S = \mathbb{R}^k$, where $S $ is the stratum of
$\mathcal{B}$ containing $\gamma(0)$.
\end{definition}
Analogous to the map \eqref{definition of tilde f} we define
\begin{equation*}
F: {T_{q_{e}}Q\times{\mathbb{R}^{k}}} \ni (v_{q_{e}},t)
\mapsto{\sum_{i=1}^{k}t_{i}F_{i}(v_{q_{e}})}\in{T_{q_{e}}Q}.
\end{equation*}
This map induces the following sets parameterized by the isotropy
type subgroups of the $G_{q_{e}}$-action on $T_{q_{e}}Q$:
\begin{equation*}
\mathcal{E}_{(H)} = \{(v_{q_e},t) \in (T_{q_e}Q)_{(H)} \times
\mathbb{R}^{k} \mid F(v_{q_e},t) \in \mathfrak{g}_{q_e} \cdot
v_{q_e} \}
\end{equation*}
\begin{equation*}
\mathcal{A}_{(H)} = \overline{\mathcal{E}_{(H)}} \cap\left(\{0\}
\times \mathbb{R}^{k}\right).
\end{equation*}
The sets $\mathcal{E}_{(H)}$ and $\mathcal{A}_{(H)} $ are
semialgebraic sets in $T_{q_e}Q \times \mathbb{R}^{k}$ and
$\mathbb{R}^k$ respectively.
\medskip
The next proposition gives an relation between the pair of subsets
$\widetilde{\mathcal{E}}_{(H)} \subset (T_{q_e}Q)/G_{q_e} \times
\mathbb{R}^k$, $\widetilde{\mathcal{A}}_{(H)} \subset
\mathbb{R}^k$  and $\mathcal{E}_{(H)} \subset (T_{q_e}Q)_{(H)}
\times \mathbb{R}^k$, $\mathcal{A}_{(H)} \subset \mathbb{R}^k$.
\begin{proposition}
\label{prop: equality between sets} Let
$\Pi:T_{q_{e}}Q\rightarrow(T_{q_{e}}Q)/G_{q_{e}}$ be the
projection onto the orbit space and $H$ and isotropy type of the
$G_{q_{e}}$-action on $T_{q_{e}}Q$. Then:
\begin{equation*}
\mathcal{E}_{(H)}=\left({\Pi}\times \operatorname{id}
\right)^{-1}\left(\widetilde{\mathcal{E}}_{(H)}\right),
\end{equation*}
\begin{equation*}
\mathcal{A}_{(H)}=\widetilde{\mathcal{A}}_{(H)}.
\end{equation*}
\end{proposition}

\noindent\textbf{Proof.} The proof is as in \cite{k}. To show the
first equality let $v_{q_e} \in \Pi^{-1}(\theta)$. Then
\begin{equation*}
\widetilde{F}(\theta,t)=0\Leftrightarrow\langle\!\langle
F(v_{q_e},t),\operatorname{grad}_{v_{q_e}}\theta_{j}(v_{q_e})\rangle\!\rangle=0,
\end{equation*}
for any $j\in\{1,\dots,l\}$. As $F(v_{q_e},t)\in
T_{v_{q_e}}(T_{q_e}Q)_{(G_{v_{q_e}})}=\mathfrak{g}_{q_e}\cdot{v_{q_e}}
\oplus{N_{v_{q_e}}^{(0)}}$ and $N_{v_{q_e}}^{(0)}$ is generated by
$\{\operatorname{grad}_{v_{q_e}}\theta_{j}(v_{q_e}),\, \text{for}
\; j\in\{1,\dots,l\}\}$, we must have
$F(v_{q_e},t)\in\mathfrak{g}_{q_e}\cdot{v_{q_e}}$. Thus
$(\theta,t)\in{\widetilde{\mathcal{E}}_{(H)}}
\Longleftrightarrow(v_{q_e},t)\in {\mathcal{E}_{(H)}}$. To prove
the second equality let $t\in{\mathcal{A}_{(H)}}$. Then there
exists $(v_{n},t_n)\in{\mathcal{E}_{(H)}}$ such that
$v_n\rightarrow{0}$ and $t_n\rightarrow{t}$. By the first
equality, this implies that
$(\Pi(v_n),t_n)\in{\widetilde{\mathcal{E}}_{(H)}}$ and the
continuity of $\Pi$ implies that
$t\in\widetilde{\mathcal{A}}_{(H)}$. Conversely, let
$t\in\widetilde{\mathcal{A}}_{(H)}$. Then there exists a sequence
$(\theta_n,t_n)\in{\widetilde{\mathcal{E}}_{(H)}}$ whose limit is
$(0,t)$. As $\mathcal{E}_{(H)}=\left({\Pi}\times \operatorname{id}
\right)^{-1}\left(\widetilde{\mathcal{E}}_{(H)}\right)$, there is
a sequence $(v_n,t_n)\in{\mathcal{E}_{(H)}}$ convergent to a point
in $(\Pi^{-1}(0),t)$. Because the action of $G_{q_e}$ on
$T_{q_e}Q$ is linear we have that $\Pi^{-1}(0)=0$ and hence
$t\in{\mathcal{A}_{(H)}}$. \quad $\blacksquare$
\medskip

\begin{remark}
If $H$ is conjugated to a proper subgroup of $G_{q_e}$ we have
$\operatorname{codim}_{\mathbb{R}^{k}}\mathcal{A}_{(H)}\geq{1}$.
\end{remark}

The next theorem, due to Koenig and Chossat \cite{kc}, gives a
sufficient criterion for the existence of symmetry breaking
branches of relative equilibria for the family of vector fields
$\overline{X_{\lambda}}$ in a neighborhood of $0_{q_{e}}$.

\begin{theorem}
\label{teorema K} Let $H$ be an isotropy group conjugate to a
proper subgroup of $G_{q_e}$ and $\overline{X_{\lambda}}\in
\mathcal{C}^{\infty}(T_{q_e}Q \times \mathbb{R},
T_{q_{e}}Q)^{G_{q_e}}_o$  verifying the transversality condition
given in  Definition \ref{transversality definition}. For
$v_{q_{e}}\in{(T_{q_{e}}Q)_{(H)}}$ define
$n_{H}:=\operatorname{rank}{N(G_{v_{q_{e}}})/G_{v_{q_{e}}}}$.
\begin{enumerate}
\item[{\rm (a)}] If
$\operatorname{codim}_{\mathbb{R}^{k}}\mathcal{A}_{(H)}\geq{2}$
then the system
$\dot{v}_{q_{e}}=\overline{X_{\lambda}}(v_{q_{e}})$ does not have
relative equilibria with isotropy conjugate to  $H$ in a
neighborhood of the origin $0_{q_{e}}$. \item[{\rm (b)}] If
$\operatorname{codim}_{\mathbb{R}^{k}}\mathcal{A}_{(H)}=1$ and
$\gamma(0)\in{\mathcal{A}_{(H)}}$ then there exists a branch of
tori of dimension $n_{H}$ bifurcating from $0_{q_{e}}$. The torus
corresponding to $\lambda$ is invariant under the flow of
$\overline{X_{\lambda}}$. The points on all of these tori have
isotropy type $H$.
\end{enumerate}
\end{theorem}

For the proof of this theorem we will need the following result of
Field and Richardson \cite{fr}:
\begin{theorem}
\label{teorema F} Let $H$ be an isotropy group conjugate to a
proper subgroup of $G_{q_e}$ and
$\widetilde{\overline{X_{\lambda}}}$ the projection on the orbit
space $(T_{q_e}Q)/G_{q_e}$ of the family of vector fields
$\overline{X_{\lambda}}\in \mathcal{C}^{\infty}(T_{q_e}Q \times
\mathbb{R}, T_{q_{e}}Q)^{G_{q_e}}_o$  verifying the transversality
condition given in  Definition \ref{transversality definition}.
\begin{enumerate}
\item[{\rm (a)}] If
$\operatorname{codim}_{\mathbb{R}^{k}}\widetilde{\mathcal{A}}_{(H)}\geq{2}$
then $\widetilde{\overline{X_{\lambda}}}^{-1}([0])$ does not
contain points with isotropy conjugate to  $H$. \item[{\rm (b)}]
If
$\operatorname{codim}_{\mathbb{R}^{k}}\widetilde{\mathcal{A}}_{(H)}=1$
and $\gamma(0)\in{\widetilde{\mathcal{A}}_{(H)}}$ then there exist
a $\mathcal{C}^{1}$ curve
$c:[0,1]\mapsto{\widetilde{\overline{X_{\lambda}}}^{-1}([0])}$
such that $c(0)=0$, $c^\prime(0)\neq{0}$, and all the points on
the curve have isotropy conjugate to $H$.
\end{enumerate}
\end{theorem}
\begin{remark}
\label{conditia fr} \normalfont In the above theorem the condition
$\overline{X_{\lambda}}\in \mathcal{C}^{\infty}(T_{q_e}Q \times
\mathbb{R}, T_{q_{e}}Q)^{G_{q_e}}_o$ is used to insure that
$f_1(0,0)=0$ and
$\left.\frac{\partial}{\partial{\lambda}}f_{1}(0,\lambda)\right|_{\lambda=0}\neq{0}$,
where $f_1$ is the first coefficient in the expansions
\eqref{x_bar} or \eqref{theta_dot_eq}. In both cases considered
here, this is equivalent to $\sigma_{X}(0)=0$ and
$\sigma_{X}^{\prime}(0)\neq{0}$.
\end{remark}
 Now we  give the proof of  Theorem \ref{teorema K}.
\medskip

\noindent\textbf{Proof.} (a)
 Suppose that the system
$\dot{v}_{q_{e}}=\overline{X_{\lambda}}(v_{q_{e}})$ admits
relative equilibria with isotropy subgroup conjugated to $H$ in a
neighborhood of $0_{q_e}$. This means that the projection of this
vector field has zeroes whose isotropy subgroups are conjugate to
$H$. By Theorem \ref{teorema F} this implies that
$\operatorname{codim}_{\mathbb{R}^{k}}\widetilde{\mathcal{A}}_{(H)}\leq
1$ and $\gamma(0)\in{\widetilde{\mathcal{A}}_{(H)}}$. This is in
contradiction with
$\operatorname{codim}_{\mathbb{R}^{k}}\widetilde{\mathcal{A}}_{(H)}=
\operatorname{codim}_{\mathbb{R}^{k}}{\mathcal{A}}_{(H)}\geq 2$.
(b) Suppose that
$\operatorname{codim}_{\mathbb{R}^{k}}{\mathcal{A}}_{(H)}=1$ and
$\gamma(0)\in{\mathcal{A}}_{(H)}$. By Proposition \ref{prop:
equality between sets}, this implies that
$\operatorname{codim}_{\mathbb{R}^{k}}\widetilde{\mathcal{A}}_{(H)}=1$
and $\gamma(0)\in{\widetilde{\mathcal{A}}_{(H)}}$. Theorem
\ref{teorema F} guarantees the existence of a $\mathcal{C}^{1}$
curve
$c:\lambda\in{[0,1]}\mapsto{\widetilde{\overline{X_{\lambda}}}^{-1}([0])}$,
$c(0)=0$, $c^\prime(0)\neq{0}$, all of whose points  have isotropy
subgroup conjugate to $H$. The inverse image of this curve under
the orbit space projection $\Pi$ is the $G_{q_e}$-orbit of a
branch of relative equilibria $v_{\lambda}$. For a fixed
$\lambda$, using Krupa's Theorem \cite{kr}, we decompose locally
in a tubular neighborhood of the orbit
$G_{q_e}\cdot{v_{\lambda}}$, the $G_{q_e}$-equivariant vector
field $\overline{X_{\lambda}}$ into its tangential and normal
components. Since the projection of this vector field to the orbit
space has $c(\lambda)$ as equilibrium we have that
$\overline{X_{\lambda}}(v_{\lambda})\in{\mathfrak{g}_{q_e}
\cdot{v_{\lambda}}}$. This means that the local normal component
of the vector field vanishes. So, the solution of the Cauchy
problem $\dot{w}_{q_{e}}=\overline{X_{\lambda}}(w_{q_{e}})$,
$w_{q_e}(0)=v_{\lambda}$, is given by
$w_{q_{e}}(t)=g(t)\cdot{v_{\lambda}}$, where
$g(t)=\operatorname{exp}(\xi_{\lambda}t)$, for small enough
$t\geq{0}$ and $\xi_\lambda\in{\mathfrak{g}_{q_e}}$. The existence
for each $\lambda$ of a flow invariant maximal torus of rank
$n_{H}$ on which the solution
$\operatorname{exp}(\xi_{\lambda}t)\cdot{v_{\lambda}}$ is
generically dense is guaranteed by Theorem 4.1 in \cite{kr}. Hence
we have obtained a branch
$w_{q_e}^{\lambda}(t)=\operatorname{exp}(\xi_{\lambda}t)\cdot{v_{\lambda}}$
of periodic or quasi-periodic solutions all of whose points have
isotropy type conjugate to $H$. \quad $\blacksquare$

Using Proposition \ref{correspondence statement}, we can deduce
the same result as the one in Theorem \ref{teorema K} for a
one-parameter family of vector fields $X_{\lambda}$ defined in a
small neighborhood around the point $q_{e}\in{Q}$. The resulting
bifurcation result will give a family of tori bifurcating from
$q_e$ whose symmetry is broken to that of the a priori given
subgroup $H \subsetneq G_{q_e}$. The bifurcating solutions are
periodic or quasi-periodic curves lying on these tori. The method
given by Proposition \ref{correspondence statement} allows one to
always translate bifurcation results obtained in the  linear
setting to analogous theorems on the manifolds itself.

In the case of the Hopf bifurcation, several additional remarks
are in order. First, the symmetry group of the problem can be
considered  to be of the form $G \times S ^1 $. This follows from
a normal form result (see e.g. \cite{gss}, Chapter XVI) which
ensures that there is a natural $S ^1$-action on $T_{q_e}Q$ such
that the family of vector fields $\overline{X_\lambda}$ are $G
\times S^1$-equivariant. Second, given $k \in \mathbb{N}$ and a
$G$-equivariant family of vector fields $\overline{X_\lambda}$ on
$T_{q_e} Q$, there is a $\mathcal{C}^k$ $\lambda$-dependent change
of variables on $T_{q_e}Q$ that maps $\overline{X_\lambda}$ into
new family of vector fields $\widetilde{X_\lambda}$ which is $G
\times S^1$-equivariant and that coincides with
$\overline{X_\lambda}$ up to order $k $. Third, one can interpret
this result as saying that close enough to $0_{q_e} $ the $S
^1$-relative equilibria of the $G \times S ^1$- equivariant family
$\widetilde{X_\lambda}$ are mapped by this change  of variables to
the set of periodic solutions of $\overline{X_\lambda}$.
Similarly, the $G \times S^1$-relative equilibria of
$\widetilde{X_\lambda}$ are mapped by this change of variables to
$G$-relative periodic orbits of $\overline{X_\lambda}$.

\section{Hamiltonian Steady State and Hamiltonian Hopf Bifurcation
with Symmetry} In this section we study a bifurcation problem
similar to the one in Theorem \ref{teorema K} in the context of a
Hamiltonian system on a symplectic vector space. To obtain the
symmetry breaking result the Hamiltonian nature of the vector
field plays a crucial role. One cannot deduce Theorem \ref{HH}
below from the generic results in \cite{kc} because the
non-degeneracy conditions are not satisfied.

Let $(V,\omega)$ be a symplectic vector space on which the compact
Lie group $G$ acts linearly and symplectically. In this case it is
known that the map $G \times V \rightarrow V $ is smooth. Let
$F_{\lambda}:V\rightarrow{\mathbb{R}}$, $\lambda\in{\mathbb{R}}$,
be a family of smooth $G$-invariant functions. If $0\in{V}$ is an
equilibrium (with symmetry $G$) of the Hamiltonian vector field
$X_{F_{0}}$, we are interested in finding conditions that
guarantee the existence of smooth branches of tori bifurcating
from zero such that all the points on these tori have symmetry
conjugate to an a priori fixed proper symmetry subgroup $H$.
Analogous to the cases considered in the previous section we shall
analyze here the following situations:

\begin{itemize}
\item \textit{Hamiltonian steady state bifurcation with symmetry}:
\begin{equation*}
D{X_{F_{\lambda}}}(0)={\sigma}(\lambda)\mathbb{J}
\end{equation*}
where $\sigma:(-\epsilon, \epsilon) \rightarrow \mathbb{R}$ is a
smooth map satisfying  ${\sigma}(0)=0$ and $\sigma^{\prime}(0)\neq
0$. In addition, it is assumed that the vector space $V $ carries
a natural complex structure $\mathbb{J}$ and that the compact Lie
group $G$ acts on this complex vector space $V$ irreducibly (as a
complex representation). In this case a basis for
$\mathfrak{gl}(V)^{G}$, the endomorphisms of $V$ which commute
with the $G$-action, is formed by $\{\operatorname{I},
\mathbb{J}\}$. This implies that a symplectic matrix in
$\mathfrak{gl}(V)^{G}$ should be a real multiple of $\mathbb{J}$.
In this case, we have passing of eigenvalues of
$D{X_{F_{\lambda}}}(0) $ on the imaginary axis at the origin as
$\lambda $ crosses zero.

\item \textit{Hamiltonian Hopf bifurcation with symmetry}:
\begin{equation*}
D{X_{F_{\lambda}}}(0)={\sigma}(\lambda)\operatorname{A_1}+{\rho}(\lambda)\operatorname{A_2}
+{\tau}(\lambda)\operatorname{A_3}+{\psi}(\lambda)\operatorname{A_4}
\end{equation*}
where
\[
\begin{array}{lll}
&\operatorname{A_1}=\left[
\begin{array}{cc}
\mathbf{0} & \mathbf{0} \\
\operatorname{I} & \mathbf{0}
\end{array}
\right], \quad &\operatorname{A_2}=\left[
\begin{array}{cc}
\mathbf{0} & \operatorname{I} \\
\mathbf{0} & \mathbf{0}
\end{array}
\right], \\ \\
&\operatorname{A_3}=\left[
\begin{array}{cc}
\operatorname{I} & \mathbf{0} \\
\mathbf{0} & -\operatorname{I}
\end{array}
\right], \quad &\operatorname{A_4}=\left[
\begin{array}{cc}
\mathbb{J} & \mathbf{0} \\
\mathbf{0} & \mathbb{J}
\end{array}
\right].
\end{array}
\]
and $\sigma, \rho, \tau, \psi :(-\epsilon, \epsilon) \rightarrow
\mathbb{R}$ are smooth maps satisfying ${\sigma}(0)=0$,
$\sigma^{\prime}(0)\neq 0$, ${\rho}(0)=-1$, $\tau(0)=0$ and
$\psi(0)=-\nu_0$, with $\nu_{0} > {0}$. In addition, it is assumed
that the vector space $V=V_{0}\oplus{V_1} $, where $V_0$ and $V_1$
are isomorphic complex dual irreducible representations of the
compact Lie group $G$. The compact Lie group $G$ is considered to
be of the form $G=\Gamma\times{S^1}$. See for example \cite{cor}.
\end{itemize}

In the following will denote by $\mathcal{C}^{\infty}(V\times
\mathbb{R},V)^{G}_{oo}$ the class of $G$-equivariant Hamiltonian
vector fields that verifies one of the above non-degeneracy
conditions.

As opposed to the general case, if the vector field is
Hamiltonian, it turns out that one needs to use only invariant
polynomials, as the proposition below shows.
\begin{proposition}
Let $(V,\omega)$ be a symplectic representation space of the
compact Lie group $G$. Denote by $\{\theta_1, \theta_2,\dots,
\theta_l\}$ a set of generators of $\mathcal{P}_G(V)$. If
$F:V\rightarrow{\mathbb{R}}$ is a smooth $G$-invariant function
then the associated Hamiltonian vector field $X_F$ has the
following expression in canonical coordinates, in terms of
invariant polynomials and their gradients:
\begin{equation}
X_F(x)=\sum_{i=1}^{l}\frac{\partial{\tilde{F}}}{\partial{\theta_i}}
(\theta_1(x),\dots,\theta_l(x))
\operatorname{\mathbb{J}}\operatorname{grad}_{x}(\theta_i)(x),
\quad\text{for any}\quad x\in{V},
\end{equation}
where $\tilde{F}:\mathbb{R}^{l}\rightarrow\mathbb{R}$ is given by:
\begin{equation*}
F(x)=\tilde{F}(\theta_1(x), \dots,\theta_l(x)),\quad\text{for
any}\quad x\in{V}.
\end{equation*}
\end{proposition}

\noindent\textbf{Proof.} The proof follows directly from the
definition of the Hamiltonian vector field. \quad $\blacksquare$
\medskip

The next step gives a method for projecting a $G$-equivariant
Hamiltonian vector field on the orbit space. To do this, let $F$
be a $G$-invariant smooth function and $X_F$ its associated
Hamiltonian vector field. Therefore, $X_F$ is $G$-equivariant on
$(V,\omega)$. Note that $\mathbb{J} \operatorname{grad} \theta_i$
are equivariant polynomials. Thus, by the proposition above we can
write $X_F$ in terms of invariants and equivariants as
\begin{equation*}
X_F(x)=\sum_{i=1}^{l}\frac{\partial{\tilde{F}}}{\partial{\theta_i}}
(\theta_1(x),\dots,\theta_l(x))\operatorname{\mathbb{J}}\operatorname{grad}_
{x}\theta_i(x),\quad\text{for any}\quad x\in{V},
\end{equation*}
where the function $\tilde{F}$ is smooth on $\mathbb{R}^{l}$. To
obtain the projection of this vector field on the orbit space
$V/G$, one has to compute the derivatives of the invariant
generators $\theta_{j}$ with respect to time:
\begin{align*}
\dot{\theta}_{j}&=\operatorname{\bold{d}}\theta_{j}(x)(\dot{x})
=(\operatorname{grad}_{x}(\theta_{j}))^{t}\cdot\left(\sum_{i=1}^{l}
\frac{\partial{\tilde{F}}}{\partial{\theta_i}}
\operatorname{\mathbb{J}}\operatorname{grad}_{x}(\theta_i)\right)
\\& =\sum_{i=1}^{l}\frac{\partial{\tilde{F}}}{\partial{\theta_i}}
\{\theta_{j},\theta_i\}.
\end{align*}
 The projection $\widetilde{X}_F$ of $X_F$ on the
orbit space $V/G$ has thus the expression
\begin{equation}
\label{hamiltonian expresion} \widetilde{X}_F(\theta)
=\mathcal{P}(\theta)\operatorname{grad}_{\theta}(\tilde{F}),
\end{equation}
where
\[
\mathcal{P}(\theta):= \left[
\begin{array}{cccc}
0 & \{\theta_1,\theta_2\} & \cdots & \{\theta_1,\theta_l\} \\
\cdots & \cdots & \cdots & \cdots \\
\{\theta_l,\theta_1\} & \cdots & \cdots & 0
\end{array}
\right].
\]
As in the Riemannian case, one search for symmetry breaking zeroes
of the projected vector field. To do this, we need a few
preliminary results. Define the following maps associated to the
family of smooth $G$-invariant functions $F:V\times{\mathbb{R}}
\rightarrow{\mathbb{R}}$
\begin{equation}
\label{definition of tilde f2} \widetilde{\Phi}:V/G \times
\mathbb{R}^{l} \ni (\theta,t) \mapsto \sum_{i=1}^{l}t_{i} \left[
\begin{array}{c}
\{\theta_1,\theta_i\} \\
\{\theta_2,\theta_i\} \\
\cdots \\
\{\theta_l,\theta_i\}
\end{array}
\right] \in V/G \subset \mathbb{R}^l
\end{equation}
\begin{equation*}
\gamma:{\mathbb{R}} \ni \lambda
\mapsto{\left(\frac{\partial{\tilde{F}}}{\partial{\theta_1}}(0,\lambda),
...,\frac{\partial{\tilde{F}}}{\partial{\theta_l}}(0,\lambda)\right)}
\in{\mathbb{R}^{l}}
\end{equation*}
and the sets
\begin{equation*}
\widetilde{\mathcal{E}}=\left\{(\theta,t)\in
V/G\times{\mathbb{R}^{l}}\mid \widetilde{\Phi}(\theta,t)=0
\right\}
\end{equation*}
\begin{equation*}
\widetilde{\mathcal{E}}_{(H)}=
\left\{(\theta,t)\in\widetilde{\mathcal{E}}\mid
{\theta\in\Pi\left(V_{(H)}\right)}\right\}
\end{equation*}
\begin{equation*}
\widetilde{\mathcal{A}}_{(H)}
=\overline{\widetilde{\mathcal{E}}_{(H)}}\cap\left(\{[0]\}\times
\mathbb{R}^{l}\right),
\end{equation*}
where $\Pi:V\rightarrow V/G$ is the orbit space projection and
$H\subset G$ is an isotropy subgroup of the $G$-action on $V$.
Recall that $\widetilde{\mathcal{E}}$ is a semialgebraic subset of
$\mathbb{R}^{2l}$. Denote by $\Sigma$ the canonical Whitney
stratification of $\widetilde{\mathcal{E}}$. This stratification
induces a Whitney semialgebraic stratification $\Sigma_{(H)}$ on
each $\widetilde{\mathcal{E}}_{(H)}$. In particular, $\Sigma_G$
stratifies
$\widetilde{\mathcal{E}}_{(G)}=\{[0]\}\times{\mathbb{R}^{l}}$.
Next we recall the notion of $G$-transversality introduced by K\oe
nig and Chossat \cite{kc} that was in turn inspired by Bierstone
\cite{bierstone} and Field \cite{field}; see Definition
\ref{transversality definition}. In the present case, this says
that the projection to the orbit space $V/G$ of a $G$-equivariant
smooth family of vector fields $\widetilde{X}_{F_{\lambda}} \in
\mathfrak{X}(V)$ is {\bfi transverse\/} to $[0]\in{V/G}$ at
$\lambda= 0\in{\mathbb{R}}$ if $\gamma: \mathbb{R} \rightarrow
\mathbb{R}^l$ is transverse at $\lambda= 0 $ to the stratum of
$\Sigma_G$ containing $\gamma(0) $, that is, $\operatorname{im}T_0
\gamma + T_{\gamma(0)} S = \mathbb{R}^l$, where $S $ is the
stratum of $\Sigma_G$ containing $\gamma(0)$. Analogous to the map
\eqref{definition of tilde f2} we define
\begin{equation*}
\Phi: V\times{\mathbb{R}^{l}} \ni (x,t)
\mapsto{\sum_{i=1}^{l}t_{i}\operatorname{\mathbb{J}}
\operatorname{grad}_{x}\theta_i(x)}\in{V}.
\end{equation*}
This map induces the following sets parameterized by the isotropy
type subgroups of the $G$-action on $V$:
\begin{equation*}
\mathcal{E}_{(H)} = \{(x,t) \in V_{(H)} \times \mathbb{R}^{l} \mid
\Phi(x,t) \in \mathfrak{g} \cdot x \},
\end{equation*}
\begin{equation*}
\mathcal{A}_{(H)} =
\overline{\mathcal{E}_{(H)}}\cap\left(\{0\}\times
\mathbb{R}^{l}\right).
\end{equation*}
 The sets $\mathcal{E}_{(H)}$ and $\mathcal{A}_{(H)} $ are
semialgebraic sets in $V \times \mathbb{R}^{l}$ and
$\mathbb{R}^l$, respectively. A proof analogous to that of
Proposition \ref{prop: equality between sets} shows that
$\mathcal{E}_{(H)}=\left({\Pi} \times \operatorname{id}
\right)^{-1} \left(\widetilde{\mathcal{E}}_{(H)}\right)$ and
$\widetilde{\mathcal{A}}_{(H)}=\mathcal{A}_{(H)}$.

\medskip

We shall prove below that if $H \neq G$, then
$\operatorname{codim}_{\mathbb{R}^{l}}\mathcal{A}_{(H)}\geq 1$.

\begin{theorem}
Let $V$ be a symplectic representation space of a compact Lie
group $G$ and $\{\theta_1,\dots,\theta_l\}$ a Hilbert base.
Suppose there exist $i_0,i_1\in\{1,\dots,l\}$, $i_0\neq{i_1}$,
$D_{i_0,i_1}:=\{x\in{V}\mid\{\theta_{i_0},\theta_{i_1}\}(x)\neq{0}\}
\neq{\varnothing}$, and $x_0 \in{V}$ such that
$\underset{x\rightarrow{x_0}}{\lim}\frac{\{\theta_{i_0},
\theta_{i}\}(x)}{\{\theta_ {i_0},\theta_{i_1}\}(x)}=0$, for any
$i\in{\{1,\dots,l\}}\setminus\{i_1\}$. Let $H$ be an isotropy
subgroup of the $G$-representation on $V$. If $(H)\neq(G)$ then we
have
$\widetilde{\mathcal{A}}_{(H)}\subset(\mathbb{R}^{1})^{\perp}=
\left\{t\in\mathbb{R}^{l}\mid t_{i_1}=0\right\}$ and
$\widetilde{\mathcal{A}}_{(G)}=\mathbb{R}^{l}$. In particular,
$\operatorname{codim}_{\mathbb{R}^{l}}
\widetilde{\mathcal{A}}_{(H)}\geq{1}$.
\end{theorem}

\noindent\textbf{Proof.} From \eqref{definition of tilde f2}, it
follows that the $i_0$-component of the map $\widetilde{\Phi}$ is
given by
\begin{equation*}
t_1\{\theta_{i_0},\theta_1\}+t_2\{\theta_{i_0},\theta_2\}
+\cdots+t_l\{\theta_{i_0},\theta_l\}.
\end{equation*}
Let $\theta: = (\theta_1, \dots, \theta_l) \in V/G \subset
\mathbb{R}^l$. If $(\theta,t)\in\widetilde{\mathcal{E}}_{(H)}$,
the above relation becomes
\begin{equation*}
t_1\{\theta_{i_0},\theta_1\}+t_2\{\theta_{i_0},\theta_2\}
+\cdots+t_l\{\theta_{i_0},\theta_l\}=0.
\end{equation*}
This implies that
\begin{equation*}
t_1\{\theta_{i_0},\theta_1\}(x)+t_2\{\theta_{i_0},\theta_2\}(x)
+\cdots+t_l\{\theta_{i_0},\theta_l\}(x)=0,\quad(\forall)
{x}\in{D_{i_0,i_1}}.
\end{equation*}
Dividing this relation by $\{\theta_{i_0},\theta_{i_1}\}(x)$ and
letting $x\rightarrow{x_0}$ one concludes that $t_{i_1}=0$. This
proves that $\widetilde{\mathcal{A}}_{(H)}\subset
\left\{t\in\mathbb{R}^{l}\mid t_{i_1}=0\right\}$. The identity
$\widetilde{\mathcal{A}}_{(G)}=\mathbb{R}^{l}$ follows directly
from the definitions. \quad $\blacksquare$

Adapting the Proposition $5.16$ of Field and Richardson \cite{fr}
to the family of vector fields $X_{F_ \lambda} $, yields the
following result. Note that $X_{F_{\lambda}}\in
\mathcal{C}^{\infty}(V \times \mathbb{R}, V)^{G}_{oo}$ implies the
condition from the Remark \ref{conditia fr}, this condition being
equivalent to $\sigma(0)=0$ and $\sigma^{\prime}(0)\neq{0}$.

\begin{theorem}
\label{HH} Let $V$ be a symplectic representation space of a
compact Lie group $G$ and $\{\theta_1,\dots,\theta_l\}$ a Hilbert
base. Suppose there exist $i_0,i_1\in\{1,\dots,l\}$,
$i_0\neq{i_1}$,
$D_{i_0,i_1}:=\{x\in{V}\mid\{\theta_{i_0},\theta_{i_1}\}(x)\neq{0}\}
\neq{\varnothing}$, and $x_0 \in{V}$ such that
$\underset{x\rightarrow{x_0}}{\lim}\frac{\{\theta_{i_0},\theta_{i}\}(x)}
{\{\theta_{i_0},\theta_{i_1}\}(x)}=0$, for any
$i\in{\{1,\dots,l\}}\setminus\{i_1\}$. Let $H$ be a given proper
isotropy subgroup of $G$ and assume that the family of Hamiltonian
vector fields $X_{F_{\lambda}}\in \mathcal{C}^{\infty}(V \times
\mathbb{R}, V)^{G}_{oo}$  satisfies the transversality condition
given in Definition \ref{transversality definition}. For
$x\in{V_{(H)}}$ define
$n_{H}:=\operatorname{rank}{N(G_{x})/G_{x}}$.
\begin{enumerate}
\item[{\rm (a)}] If
$\operatorname{codim}_{\mathbb{R}^{l}}\mathcal{A}_{(H)}\geq{2}$
then the system $\dot{x}=X_{F_{\lambda}}(x)$ does not have
relative equilibria with isotropy conjugate to  $H$ in a
neighborhood of the origin $0$. \item[{\rm (b)}] If
$\operatorname{codim}_{\mathbb{R}^{l}}\mathcal{A}_{(H)}=1$ and
$\gamma(0)\in{\mathcal{A}_{(H)}}$ then there exists a branch of
tori of dimension $n_{H}$, formed by points of isotropy type $H$,
bifurcating from $0$. The torus corresponding to $\lambda$ is
invariant under the flow of $X_{F_{\lambda}}$.
\end{enumerate}
\end{theorem}

\medskip

We shall make some remarks regarding the Hamiltonian Hopf
bifurcation. We have made the hypothesis that the symmetry group
of the problem is of the form $G \times S^1$. This always holds if
we assume only that the system has a compact Lie group $G$ of
symmetries represented on the symplectic vector space $V$. The
reason for this is that in the case of the Hamiltonian Hopf
bifurcation there is a natural $S^1$-action on $V$ given by the
exponential of the semisimple part of the linearization on the
resonance space corresponding to $\nu_0$. More precisely, let $(V,
\omega, F_\lambda)$ be a $\lambda$-parameter family ($\lambda \in
\mathbb{R}$ in a neighborhood of the origin) of $G$-Hamiltonian
systems such that $F_0 (0)=0$, $\mathbf{d} F_0 (0)=0$, and the
$G$-equivariant infinitesimally symplectic linear map
$DX_{F_0}(0)$ is nonsingular and has $\pm i\nu_0$ as eigenvalues.
Let $(U_{\nu_0}, \omega|_{U_{\nu_0}})$ be the resonance space of
$DX_{F_0}(0)$ with primitive period $T_{ \nu_0} := 2 \pi/ \nu_0$.
Recall that the resonance space  is the direct sum of the real
eigenspaces of $DX_{F_0}(0)$ corresponding to eigenvalues of the
form $\pm k i \nu_0 $ with $k = 1, 2, \dots$. A theorem of
Vanderbauwhede and van der Meer \cite{vvdm} states that for each
$k\geq 0$ one can find a $C^k$-mapping $\psi:U_{\nu_0}\times
\mathbb{R}\rightarrow V$ and a $C^{k+1}$-function
$\widehat{F_\lambda}:U_{\nu_0}\times
\mathbb{R}\rightarrow\mathbb{R}$ such that  $\psi(0,\lambda)=0$,
for all $\lambda$,
$D_{U_{\nu_0}}\psi(0,0)=\mathbb{I}_{U_{\nu_0}}$,  and
$\widehat{F_\lambda}$ is a $G\times S^1$-invariant function  that
coincides with $F_\lambda$ up to order $k+1$. In the  same paper
it is shown (Theorem 3.2) that if $\lambda$ is close to $0$ then
the $S^1$-relative equilibria  of the $G\times S^1$-invariant
Hamiltonian $\widehat{F_\lambda}$ that lie in a sufficiently small
neighborhood of the origin in $U_{ \nu_0}$  are mapped by
$\psi(\cdot,\lambda)$ to the set of periodic solutions of $(V,
\omega, F_\lambda)$ in a neighborhood of $0\in V$ with periods
close to $T_{ \nu_0}$. Thus, on can replace the search of periodic
orbits of $(V, \omega, F_\lambda)$ by the search of $S^1$-relative
equilibria of the $G\times S^1$-invariant family of Hamiltonian
systems $(U_{ \nu_0}, \omega|_{U_{ \nu_0}}, \widehat{F_\lambda})$.
Similarly, the $G \times S^1$-relative equilibria of
$X_{\widehat{F_\lambda}}$ are mapped by $\psi( \cdot , \lambda) $
to $G$-relative periodic orbits of $X_{F_\lambda}$.

\section{Bifurcation of relative equilibria in Hamiltonian systems on symplectic vector spaces}
In this section we will prove a bifurcation result in symplectic
vector spaces in the context of Hamiltonian systems. We consider a
symplectic representation space $(V,\omega)$ of a compact Lie
group $G$ and let $\{\theta_1,\dots,\theta_l\}$ a Hilbert base. We
will denote by $V_{(H)}\subset{V}$ the principal stratum of the
action of $G$ on $V$ and by $\Pi:V\rightarrow{V/G}$ the orbit
space projection. Recall that each connected component of
$V_{(H)}$ is an open dense subset of the vector space $V$. We
consider a family of $G$-invariant smooth functions
$F_{\lambda}:V\rightarrow{\mathbb{R}}$, $\lambda\in\mathbb{R}$,
such that there exists $\lambda_0\in{\mathbb{R}}$ and
$v_0\in{V_{(H)}}$ a relative equilibrium for the Hamiltonian
vector field $X_{F_{\lambda_0}}$. We will prove that under a
certain non-degeneracy condition we have the existence of a smooth
branch of relative equilibria (starting from $v_0$) for the
Hamiltonian vector fields $X_{F_{\lambda}}$ for $\lambda$ close to
$\lambda_0$. To do this we recall some of the techniques necessary
to find the expression of the projected Hamiltonian vector field
on $V/G$ . More precisely, if
$F_{\lambda}(x)=\tilde{F}_{\lambda}(\theta_1(x),\dots,\theta_l(x))$,
for any $x\in{V}$, where
$\tilde{F}_{\lambda}:\mathbb{R}^{l}\rightarrow{\mathbb{R}}$ is
smooth, then the $j$-component of the projection on the orbit
space $V/G$ of the Hamiltonian vector field $X_{F_{\lambda}}$ is
given by:
\begin{equation*}
\sum_{i=1}^{l}\frac{\partial{\tilde{F}_{\lambda}}}{\partial{\theta_i}}
\{\theta_{j},\theta_i\}.
\end{equation*}
Searching for relative equilibria of a vector field is the same as
searching for equilibria for the projected vector field on the
orbit space.  Thus, we reduced the problem to find branches in
$\lambda$ of solutions of the system:
\begin{equation*}
\sum_{i=1}^{l}\frac{\partial{\tilde{F}_{\lambda}}}{\partial{\theta_i}}
(\theta_1,\dots,\theta_l)\{\theta_{j},\theta_i\}=0,
\quad\text{for}\quad{j\in\{1,\dots,l\}}.
\end{equation*}
To do this, we introduce the function
$g:(V/G\subset{\mathbb{R}^{l}})\times{\mathbb{R}}\rightarrow{\mathbb{R}^{l}}$,
$g=(g_1,\dots,g_l)$, where the components $g_j$ are defined by:
\begin{equation*}
g_j((\theta_1,\dots,\theta_l),\lambda):=\sum_{i=1}^{l}\frac{\partial{\tilde{F}_
{\lambda}}}{\partial{\theta_i}}
(\theta_1,\dots,\theta_l)\{\theta_{j},\theta_i\},
\quad\text{for}\quad{j\in\{1,\dots,l\}}.
\end{equation*}
As $v_0\in{V_{(H)}}$ is a relative equilibrium of
$X_{F_{\lambda_0}}$ this implies that $\theta_{0}=:\Pi(v_0)$ is an
equilibrium point for the projection on the orbit space
$V/G\subset{\mathbb{R}^{l}}$ of the Hamiltonian vector field
$X_{F_{\lambda_0}}$. This means that
$\theta_{0}\in{V/G}\subset{\mathbb{R}^{l}}$ is a solution of the
equation $g(\theta,\lambda_0)=0$. The next theorem will give a
sufficient condition that guarantees the existence of a smooth
curve of solutions of the equation $g(\theta,\lambda)=0$, in an
open neighborhood of the point $(\theta_{0},\lambda_{0})$.
\begin{theorem}
\label{HVS} Let $(V,\omega)$ be a symplectic representation space
of a compact Lie group $G$ and $\Pi:V\rightarrow{V/G}$ be the
orbit space projection. Let
$F_{\lambda}:V\rightarrow{\mathbb{R}}$, $\lambda\in\mathbb{R}$, be
a family of $G$-invariant smooth functions. Suppose that
$v_0\in{V}$ is a relative equilibrium of the Hamiltonian vector
field $X_{F_{\lambda_0}}$ which belongs to the principal stratum
of the $G$-action on $V$. Denote $\theta_{0}=:\Pi(v_0)$. If
$\operatorname{\bold{d}}_{\theta}g(\theta_{0},\lambda_{0})$ is
non-degenerate, then there exists an open subset $I$ of
$\mathbb{R}$ containing $\lambda_0$ and a smooth curve
$\gamma:\lambda\in{I}\mapsto
{\gamma(\lambda)}\in{V/G}\subset{\mathbb{R}^{l}}$, such that
$\gamma(\lambda)\in{V/G}\subset{\mathbb{R}^{l}}$ lifts to relative
equilibria of $X_{F_{\lambda}}$ and
$\gamma(\lambda_0)=\theta_{0}\in{V/G}\subset{\mathbb{R}^{l}}$.
\end{theorem}

\noindent\textbf{Proof.} The proof is a consequence of the
implicit function theorem and of the fact that
$v_{0}\in{V_{(H)}}\subset{V}$ is an element of the principal
stratum; the connected component containing it is an open subset
of $V$. \quad $\blacksquare$

Note that in this case the only thing we can say about the
symmetry of the relative equilibria on the branch is that for
$\lambda$ small enough all of them are of principal type. This
theorem can be used in the following  way to study bifurcating
branches of relative equilibria from a given one with isotropy
$K$. Let $v_0 \in V $ be a relative equilibrium with isotropy $K$.
Since $G$ is compact, so is $K $ and hence the subspace $V^K$ of
$K $-fixed vectors is a symplectic subspace of $(V, \omega)$. Let
us assume that the normalizer $N(K)$ of $K $ in $G $ is closed,
which implies that $N(K)/K $ is a compact Lie group that acts
naturally on $V^K $. Note that the restriction of this action on
the $K$-invariant open subset $V_K$ is free. Then we can apply the
previous theorem to $V^K $, $N(K)/K $, and the restrictions of the
vector fields $X_{F_\lambda}$ to the open subset $V_K$ of $V^K$.
The corresponding  function $g $ from before will be restricted to
the open set $V_K/(N(K)/K)$ and, under the same\ non-degeneracy
hypothesis, we obtain $\lambda$-branches of relative equilibria
for $X_{F_\lambda}$ consisting  of points with isotropy exactly
$K$.

\section{From symplectic vector spaces to manifolds}
 In this short section we will give a method of translating
the results obtained in the previous two sections from symplectic
vector spaces to symplectic manifolds. As all the symmetry
breaking results in sections $5$ and $6$ are of local type we will
make use of the Equvariant Darboux Theorem to pass locally from a
symplectic manifold to a symplectic vector space.
\begin{theorem}
\label{Darboux} Suppose that $\Gamma$ is a compact Lie group
acting on a finite-dimensional symplectic manifold $Q$ and let
$\omega$ be a $\Gamma$-invariant symplectic form on $Q$. Let
$q\in{Q}$ be a $\Gamma$-invariant point. Then there are open
$\Gamma$-invariant neighborhoods $U \subset Q$ of $q$ and $V
\subset T_q Q$ of $0_q \in T_q Q$ and a $\Gamma$-invariant
symplectic diffeomorphism $f:(U, \omega) \rightarrow (V,
\omega(q))$, such that $f(q) = 0_q $.
\end{theorem}
This theorem follows from the $G$-relative Darboux Theorem (see,
for example, \cite{or}).
\begin{theorem}
Let $Q$ be a smooth manifold and $\omega_0$ and $\omega_1$ two
symplectic forms on it. Let $G$ be a Lie group acting properly on
$Q$ and symplectically with respect to both $\omega_0$ and
$\omega_1$. Assume that
\begin{equation*}
\omega_{0}(g\cdot{q})(v_{g\cdot{q}},w_{g\cdot{q}})=
\omega_{1}(g\cdot{q})(v_{g\cdot{q}},w_{g\cdot{q}})
\end{equation*}
for all $g\in{G}$ and
$v_{g\cdot{q}},w_{g\cdot{q}}\in{T_{g\cdot{q}}Q}$. Then there exist
two open $G$-invariant neighborhoods $U_0$ and $U_1$ of
$G\cdot{q}$ and a $G$-equivariant diffeomorphism
$\Psi:U_{0}\rightarrow{U_{1}}$ such that
$\Psi\mid_{G\cdot{q}}=\operatorname{Id}$ and
$\Psi^{*}\omega_1=\omega_0$.
\end{theorem}

Theorem \ref{Darboux} is proved as follows. On the symplectic
manifold $(Q,\omega)$ consider around the point $q$ an open
neighborhood $W\subset{Q}$ contained in a Darboux chart. Thus, on
this neighborhood the symplectic form $\omega_1$ is the pull back
of the symplectic form $\omega(q)$ on $T_qQ$. Recall that
$\omega(q)$ is $\Gamma$-invariant, where $\Gamma$ acts by
linearization of the original action on $T_qQ$. Since the point
$q$ is $\Gamma$-invariant we can shrink the neighborhood $W$ to a
$\Gamma$ invariant open neighborhood $B\subset W$ of $q$. So, on
the $\Gamma$-invariant manifold $B$ we have two $\Gamma$-invariant
symplectic forms: $\omega$ and $\omega_1$. Now using the
$\Gamma$-relative Darboux theorem we get the existence of two open
$\Gamma$-invariant neighborhoods $U_{0}\subset{B}$ and
$U_{1}\subset{B}$ of $q$ and a $\Gamma$-equivariant diffeomorphism
$\Psi:U_{0}\rightarrow{U_{1}}$ such that
$\Psi\mid_{\Gamma\cdot{q}}=\operatorname{Id}$ and
$\Psi^{*}\omega=\omega_1$.

 The nonlinear version of the results from sections $5$ and $6$
can be stated in the following theorem.
\begin{theorem}
Let $G$ be a Lie group acting properly on a symplectic finite
dimensional manifold $(Q,\omega)$. Let
$F_{\lambda}:Q\rightarrow{\mathbb{R}}$, $\lambda\in{\mathbb{R}}$,
be a family of $G$-invariant smooth functions and $q_{e}\in{Q}$ be
a relative equilibrium of the Hamiltonian vector field
$X_{F_{\lambda_0}}$. Then, in the conditions of Theorems \ref{HH}
and \ref{HVS}, the same conclusions hold on the manifold $Q$
around the point $q_e$.
\end{theorem}

\noindent\textbf{Proof.} The proof is a consequence of  Theorem
\ref{Darboux} by considering the compact Lie group $G_{q_e}$
acting on the symplectic manifold $(Q,\omega)$. \quad
$\blacksquare$

\noindent {\sc P. Birtea} \\
Departamentul de Matematic\u a, Universitatea de Vest,
RO--1900 Timi\c soara, Romania.\\
Email: {\sf birtea@math.uvt.ro}
\medskip

\noindent {\sc M. Puta} \\
Departamentul de Matematic\u a, Universitatea de Vest,
RO--1900 Timi\c soara, Romania.\\
Email: \texttt{puta@math.uvt.ro}
\medskip

\noindent {\sc T.S. Ratiu}\\
Centre Bernoulli, {\'E}cole Polytechnique F{\'e}d{\'e}rale de
Lausanne, CH--1015 Lausanne,
Switzerland.\\
Email: {\sf tudor.ratiu@epfl.ch}
\medskip

\noindent {\sc R. M. Tudoran}\\
Centre Bernoulli, {\'E}cole Polytechnique F{\'e}d{\'e}rale de
Lausanne, CH--1015 Lausanne, Switzerland; Departamentul
deMatematic\u a, Universitatea de Vest,
RO--1900 Timi\c soara, Romania.\\
Email: {\sf razvan.tudoran@epfl.ch}

\end{document}